\documentclass{article}
\usepackage{arxiv}
\usepackage{xparse}
\usepackage[pdftex]{graphicx}
\usepackage{bm}
\usepackage{amsmath,amsthm,amsfonts}
\usepackage{algorithm}
\usepackage{algorithmic}

\usepackage[caption=false,font=normalsize,labelfont=sf,textfont=sf]{subfig}
\usepackage{tikz}
\usepackage{siunitx}
\usetikzlibrary{calc}
\usetikzlibrary{shapes,arrows,positioning}
\tikzstyle{decision} = [diamond, draw, fill=blue!20, 
    text width=6em, text badly centered, node distance=3cm, inner sep=0pt]
\tikzstyle{block} = [rectangle, draw, fill=blue!20, 
    text width=8em, text centered, rounded corners, minimum height=4em]
\tikzstyle{line} = [draw, -latex']
\tikzstyle{cloud} = [draw, ellipse,fill=red!20, node distance=3cm,
    minimum height=2em]
\NewDocumentCommand\norm{mg}
{\|{#1}\|\IfNoValueTF{#2}{}{_{#2}}}

\newcommand\bSigma{\mathbf{\Sigma}}
\newcommand\bx{\mathbf{x}}
\newcommand\bv{\mathbf{v}}
\newcommand\bK{\mathbf{K}}

\newcommand\bM{\mathbf{M}}
\usepackage{lineno}
\begin{document}
\title{The asymmetric particle population density method for simulation of coupled noisy oscillators}
\author{Ningyuan~Wang
        and~Daniel~B~Forger
\thanks{N. Wang is with the Department
of Mathematics, University of Michigan, 
Ann Arbor, MI, 48109}
\thanks{D. B. Forger is with Department of Mathematics, Department of Computational Medicine and Bioinformatics and Michigan Institute for Data Science, University of Michigan,
Ann Arbor, MI, 48109
 e-mail: forger@umich.edu
}
}

\markboth{ }%
{Wang, Forger: APPD for coupled noisy oscillators}

\maketitle

\begin{abstract}
A wide variety of biological phenomena can be modeled by the collective activity of a population of individual units. A common strategy for simulating such a system, the population density approach, is to take the macroscopic limit and update its population density function. However, in many cases, the coupling between the units and noise gives rise to complex behaviors challenging to existing population density approach methods. To address these challenges, we develop the asymmetric particle population density (APPD) method that efficiently and accurately simulates such populations consist of coupled elements.
The APPD is well-suited for a parallel implementation. We compare the performance of the method against direct Monte-Carlo simulation and verify its accuracy by applying it to the well-studied Hodgkin-Huxley model, with a range of challenging scenarios. We find that our method can accurately reproduce complex macroscopic behaviors such as inhibitory coupling-induced clustering and noise-induced firing while being faster than the direct simulation. 
\end{abstract}


\section{Introduction}
Many biological behaviors can be modeled by coupled oscillating elements. Examples on various scales include the rhythmic flashing of fireflies~\cite{buck1988synchronous}, neurons~\cite{hodgkin1952measurement}, and the circadian rhythm~\cite{to2007molecular}. Understanding the macroscopic behavior (e.g., synchronization of the population) for a large population without calculating the dynamics of all oscillators individually is of great value. A common approach to this problem is to consider the population density distribution~\cite{winfree1967biological, jansen1995electroencephalogram, sanger1996probability, nykamp2000population, kuramoto2003chemical}. However, most existing methods also have severe constraints on the oscillator model, such as requiring the oscillator to closely follow a limit cycle, or that the coupling is weak. These constraints and simplifications can leave some important applications behind, such as performing a parameter search to tune the parameters that generate desired macroscopic properties. Additionally, since noise can also generate macroscopic patterns~\cite{goldwyn2011and, ko2010emergence}, it needs to be accounted for accurately. The goal of our proposed method is to form a middle ground between the direct simulation of all oscillators and existing population density-based methods, such that: 1)~the full dynamics of the oscillator are retained; 2)~accurately accounts for the noise; 3)~accommodates complicated coupling between the oscillators, and; 4)~faster than the direct simulation.

Here, we present the asymmetrical particle population density (APPD) method to simulate such a population of noisy all-to-all-coupled oscillators. Our method relies on a population density approach, which tracks the probability that any element is in any particular state at a time. While such methods are widely considered computationally prohibitive for biologically realistic models, Stinchcombe and Forger~\cite{STINCHCOMBE2016} note that most of the models are dissipative. Therefore, only a fraction of the possible system states has a non-negligible number of oscillators at any given time. Additionally, coupling between the oscillators can concentrate the distribution. As demonstrated in~\cite{STINCHCOMBE2016}, these distributions can be simulated using particle methods without any further mathematical simplifications. However, the particle method presented in~\cite{STINCHCOMBE2016} is still limited by two factors: 1) Particles are symmetric in space, which covers a trajectory inefficiently. 2) The computation of the noise term is based on the interaction between particles~\cite{mas2002diffusion}, thus requires more computation and is hard to parallel than if particles can be updated independently. We address these issues to create a fast and efficient method.

While APPD is inspired by the previous work~\cite{STINCHCOMBE2016} of using a particle method in a similar setting, it is novel in that the particles are asymmetric, allowing it to change its shape with the flow. This allows the population density to be very narrow in some dimensions. In practice, we find this true in neuronal simulations, which greatly reduces the number of particles needed and increases these methods' accuracy. Additionally, the particles can be updated individually except for coupling, allowing a much easier and more efficient parallel implementation of the algorithm. 

The practice of using an asymmetric deforming kernel to solve a convection-diffusion equation or a Fokker-Planck equation is not new, either in the context of convection-diffusion equation~\cite{rossi1996resurrecting,rossi2005achieving}, or in a very different context of Kalman-Bucy filtering~\cite{sarkka2007unscented,arasaratnam2010cubature,gustafsson2011some,kulikov2017accurate,wang2021level}. The concept of splitting and combining such particles is also known~\cite{xie2014adaptive,abgrall2014immersed,berchet2021adaptive}. But our main innovation is as follows: First, the derivation of our method is based on tracking movement of a level set, which is very different from existing methods based on a direct computation of terms in the covariance matrix. As a result, our method naturally tracks a~\emph{square-root} of the covariance, which gives improved numerical stability, as known in the Kalman filter community~\cite{sarkka2007unscented,arasaratnam2010cubature} but has not been introduced to this context. This improved stability is needed as previous particle methods are restricted to $2$-dimensional scenarios, whereas the Hodgkin-Huxley model~\cite{hodgkin1952measurement} has a $4$-dimensional state space. Additionally, our method does not require a Taylor expansion in time with an \emph{a priori} timestep to represent diffusion, which enables an adaptive timestep method to be directly applied. Also, the splitting of a single into three particles with optimized weight and distance introduces less error than a 2-particle split, as in \cite{berchet2021adaptive}, for example. For the sake of simplicity, in the following discussion, we will refer to elements as oscillators, a key target application. However, this method would also apply to systems that are not strictly oscillators, such as being quasiperiodic.

\section{Method}
\subsection{Overview of the method:} 
For simulating a large population of coupled noisy oscillators, we consider the direct Monte-Carlo simulation as the baseline method that we compare against. In the direct Monte-Carlo method, each oscillator is updated individually based on its dynamics, noise, and coupling. In our proposed APPD, the intuition for the method is simple: a particle not only represents oscillators exactly at this state, but also nearby oscillators with some locally linear approximation. In the following parts, we present how such intuition is implemented, namely: how the dynamics and noise of the system are updated for a single particle in a deterministic algorithm, how the local linearity is preserved by splitting the particle, and other technicalities that arise during the process.
 
In the APPD, a single particle is a weighted multivariate Gaussian. Simulating the probability distribution function (PDF) of a similar stochastic process using a single multivariate Gaussian is a well-studied topic, in the context of Kalman-Bucy filtering~\cite{kalman1961new} and its generalizations~\cite{sarkka2007unscented,arasaratnam2010cubature,gustafsson2011some,kulikov2017accurate}. These methods utilize the fact that a Gaussian function is preserved if the dynamics of the system is linear in space. Therefore, when the PDF is concentrated in a small region, it can be considered approximately locally linear, and the Gaussian is preserved. However, these existing methods~\cite{sarkka2007unscented,arasaratnam2010cubature,gustafsson2011some,kulikov2017accurate} are based on It\^{o}-Taylor expansion of the underlying equation, thus the representation of noise is dependent on the choice of the timestep size. In \cite{wang2021level}, we proposed the Level Set Kalman Filter (LSKF), whose representation is independent of the time discretization, and has superior accuracy than the Continuous-Discrete Cubature Kalman Filter.
Therefore we base the particle in our purposed method on the \emph{time-update} of the LSKF.

While the Fokker-Planck equation for the LSKF and the convection-diffusion equation for the population takes the same form, the challenges for porting LSKF, or similar multivariate-Gaussian-based methods to the context of coupled noisy oscillators are nontrivial. First, oscillators can have a limit cycle with strong contraction, resulting in the PDF being very thin in certain directions, corresponding to a close-to-singular covariance matrix in a multivariate-Gaussian approximation. Therefore the method has to be robust with respect to singular covariance.
For the LSKF, since it works with a square root of the covariance matrix and is robust for semi-positive definite covariances, it is stable as long as a backward-stable linear equation solver is used, and the near-singular eigenvectors are dropped. The next challenge is that the approximately linear condition for the dynamics is not valid for oscillators; additionally, we are interested in more complicated distributions than what can be described by a single multivariate Gaussian. This challenge is resolved by introducing a scheme to split and combine particles.

\subsection{Problem formulation}
The subject of the simulation is a population consisting of identical oscillators subject to dynamics described by a velocity field $\mathbf{v}$, and a noise that can be described by a Wiener process. 
We assume that an oscillator can be described by a $d$-dimensional state variable. Suppose that the oscillators are also subject to noise, then this probability density $u$ is evolved according to the following convection-diffusion equation:
\begin{equation}
\label{eqn:fokker_planck}
\frac{\partial u}{\partial t} =  \nabla \cdot \bK \nabla u - \nabla \cdot (\mathbf{v} u),
\end{equation}
in which $\bK$ is a constant \textbf{diffusion matrix}, and $\mathbf{v} = \mathbf{v}(\mathbf{x},u)$ is the \textbf{convection velocity} that corresponds to the ordinary differential equation of the oscillator, and the coupling effect between oscillators. 

We further assume that the population is large (in fact, infinitely large) such that individual fluctuations will not affect the population. Hence the probability distribution and the population distribution is the same. This assumption is essential as the coupling would be more complex otherwise.

To accommodate coupling, we assume the velocity field takes the following form:
\begin{equation}
\mathbf{v} = \mathbf{v}_d(\mathbf{x}) + \mathbf{v}_c(u,\mathbf{x}),
\end{equation}
in which $\bv_d$ is the \textbf{dynamic velocity} corresponding to the dynamics that are determined by the state variable $\mathbf{x}$ only. $\mathbf{v}_c$ is the \textbf{coupling velocity} corresponding to the coupling term that depends on the state variable $\mathbf{x}$, and the population distribution $u$ at the same time. Moreover, we assume the coupling velocity $\mathbf{v}_c$ is of the following form:
\begin{equation}
\label{eqn:coup_vel}
\mathbf{v}_c(u,\mathbf{x}) = \mathbf{w}_c(\mathbf{x}) \cdot \mathbf{L}[u],
\end{equation}
in which $\mathbf{w}_c(\cdot)$ is a vector-valued function, $L[\cdot]$ is a linear functional on the PDF $u$.

In principle, the dependence of coupling velocity on the density makes (\ref{eqn:fokker_planck}) a nonlinear differential equation, which gives rise to complex behaviors. However, since the coupling velocity $\mathbf{v}_c$ as defined in (\ref{eqn:coup_vel}) depends on global quantities of $u$ over the domain, and can be assumed to change rather slowly for our applications. Therefore, we approximate $\mathbf{v}_c = \mathbf{v}_c(\mathbf{x})$ to be independent of $u$ during the update, except when the velocity function $\mathbf{v}_c$ is updated at some fixed timesteps. This approximation keeps (\ref{eqn:fokker_planck}) linear in the following derivations.

To numerically simulate the population density, a discretization is required. In a \textbf{particle method}, the population density is discretized by approximating the population density as a linear combination of \textbf{particles}:
\begin{equation}
u(\mathbf{x}) \approx \sum_{i \in I} w_i K_i(\mathbf{x-x_i}),
\end{equation}
where for each particle with index $i$, $K_i(\mathbf{x})$ is the \textbf{kernel function} of the particle, $w_i$ is the \textbf{weight} of the particle, and $\mathbf{x_i}$ is the \textbf{center location} of the particle.

In our method, each particle is a weighted multivariate Gaussian. In particular, each Gaussian particle is represented by its \textbf{weight} $w_i$, \textbf{center location} $\mathbf{x}_i$ and \textbf{covariance matrix} $\bSigma_i$, and its \textbf{kernel function} $K_i$ is given by: 
\begin{equation}
K_i(\mathbf{x}) = \frac{1}{\sqrt{(2 \pi )^d\det(\bSigma_i)}}\exp\left(-\frac{1}{2} \bx^T \bSigma^{-1}_i \bx\right).
\label{eqn:kernel_function}
\end{equation}
To improve the numerical stability and derive a simpler update algorithm, we track and store a \textbf{square root} of the covariance matrix $\bSigma = \bM \bM^T$ instead, and only evaluate the covariance matrix when needed.
\subsection{Single particle update}
This section describes the particle update algorithm, which is a special case of the Level Set Kalman Filter method as derived in \cite{wang2021level}. 
Only the description is included here. Therefore readers interested in the derivation and proof should refer to \cite{wang2021level}.

Recall (\ref{eqn:fokker_planck}):
\begin{equation}
\label{eqn:heatadv}
\frac{\partial u}{\partial t} =  \nabla \cdot \bK \nabla u - \nabla \cdot (\bv u).
\end{equation}

For a single particle centered at $\bx_0$, the probability density function is given by
\begin{equation}
u(\bx,0) = \frac{1}{\sqrt{(2 \pi)^d \det(\bSigma)}} \exp\left(-\frac{(\bx-\bx_0)^T \bSigma^{-1}(\bx-\bx_0)}{2}\right).
\end{equation}
Consider the \textbf{level set} of the function $F$: 
\begin{equation}
\label{eqn:levelset}
\left\{\bx | F(\bx,t) := \frac{u(\bx,t)}{u(\mathbf{0},t)} = c \right\} \quad (0<c<1).
\end{equation}
For $u$ being a Gaussian particle, all level sets are ellipsoids. Tracking the movement of the Gaussian particle is equivalent to tracking one of its ellipsoid level sets as defined in (\ref{eqn:levelset}). A factorization of the covariance matrix $\bSigma = \bM \bM^T$ (where $\bM$ is called a \textbf{square root} of $\bSigma$) represents a level set in the sense that the column vectors lie on the same ellipsoid. We represent the columns of the matrix $\bM$ as
$\bM_i$, and advance them in time. The ordinary differential equation is given by:
\begin{equation}
\label{eqn:level_set_ODE_num}
\frac{\partial \bM_i}{\partial t} = \frac12\left(\mathbf{\bv}(\bx_0+\bM_i) - \bv(\bx_0-\bM_i)\right) + \bK(\bM^{T})^{-1}\mathbf{e}_i,
\end{equation}
where $\bx_0$ is the \textbf{center location} of the Gaussian particle, and the matrix inverse shall be understood as solve the equation with a backward-stable solver.

A matrix short-hand form is as follows (where the matrix-vector additions are defined entrywise): 
\begin{equation}
\label{eqn:level_set_ODE_num_matform}
\frac{\partial \bM}{\partial t} = \frac12\left(\mathbf{\bv}(\bx_0+\bM) - \bv(\bx_0 -\bM)\right) + \bK(\bM^{T})^{-1},
\end{equation}
and the velocity for center is given by: 
\begin{equation}
\frac{\partial \bx_0}{\partial t} = \frac{1}{2d}\sum_{i=1}^{d}{\bv(\bx_0 + \bM_i) + \bv(\bx_0 - \bM_i)}.
\end{equation}
Notice that to evaluate the velocity for one point $\bx_i$ on the level set, the center location, and all other points $\bx_j$ for this Gaussian kernel are needed, hence the points on a level set cannot be updated independently, and the dimension for the ODE solver is $d\times (d+1)$. However, the updates of different Gaussian particles are independent. 

As a summary, a particle is updated as follows: 
\begin{algorithm}[H]
\caption{Particle Update}
\label{alg:particle_update}
\begin{algorithmic}
\REQUIRE{center $\bx$, and a square root $\bM$ of the covariance matrix at previous time step.}
\STATE{ Pass $\bM$ and $\bx$ as the state variable to an ODE solver with derivative defined by  (\ref{eqn:level_set_ODE_num_matform}).}
\RETURN{center $\bx'$, and  a square root $\bM'$ of the covariance matrix.}
\end{algorithmic}
\caption{summary of particle update}
\end{algorithm}

\subsection{Splitting a particle}
Since updating the center location and the square root of the covariance matrix for a single particle uses a linear approximation of the convection velocity, we need a method to control the size of the Gaussian particle that would otherwise expand due to diffusion. Though the Gaussian particles have infinite support, most of the density is within $2$ standard deviations in each eigenvector direction, and we consider this region to be its effective support.
The criterion we use for checking linear approximation is as follows: we check the soundness of the linear approximation in the direction of the off-center points used in the particle update. We decide if the particle needs to be split if this particular relative error for any of the off-center points is larger than a threshold chosen by the user: 
\begin{equation}
\epsilon = \frac{\norm{(\bv(\bx +2 \Delta \bx) - \bv(\bx)) - 2 (\bv(\bx+ \Delta \bx -\bv(\bx)))}}{2\norm{\bv(\bx)}},
\end{equation}where $\bx$ is the center of the particle, and $\Delta \bx$ is the offset of the off-center point. For the numerical examples that follow, the value is chosen as $0.05$.

Once the velocity field $v$ in the effective support deviates from the linear approximation larger than this tolerance, we need a method to reduce the covariance, thereby reducing the effective support of the particle. 
We achieve this by splitting the particle into three thinner particles in the corresponding direction with half variance in that direction. 

 Up to some rotation, the kernel function (\ref{eqn:kernel_function}) can be reformulated as: 
\begin{equation}
K_i(\mathbf{x}) = \frac{1}{(2 \pi)^{d/2} \sqrt{\sigma_1\sigma_2\dots \sigma_d}}\exp\left(-\frac{1}{2} \left(\frac{x_1^2}{\sigma_1}+\frac{x_2^2}{\sigma_2}+\dots+\frac{x_d^2}{\sigma_d}\right)\right),
\end{equation}
in which $\sigma_1 \geq \sigma_2 \geq \dots \geq \sigma_d $ are the eigenvalues of $\bSigma$. We here consider a method to reduce the width of the population in $x_1$ direction by representing the original population as a summation of three particles with variance in $x_1$ direction reduced by half. The particle in the center would weigh $1-2\omega$, and have its location the same as what it replaced, while the two off-center particles would both weigh $\omega$, and have its center location shift by $a$ and $-a$ in the direction of $x_1$ respectively. Note the weight of the three children particles add up to $1$, and the total weight of the particle is conserved after the operation.

\begin{figure}[H]
\centering
\includegraphics[width =\textwidth]{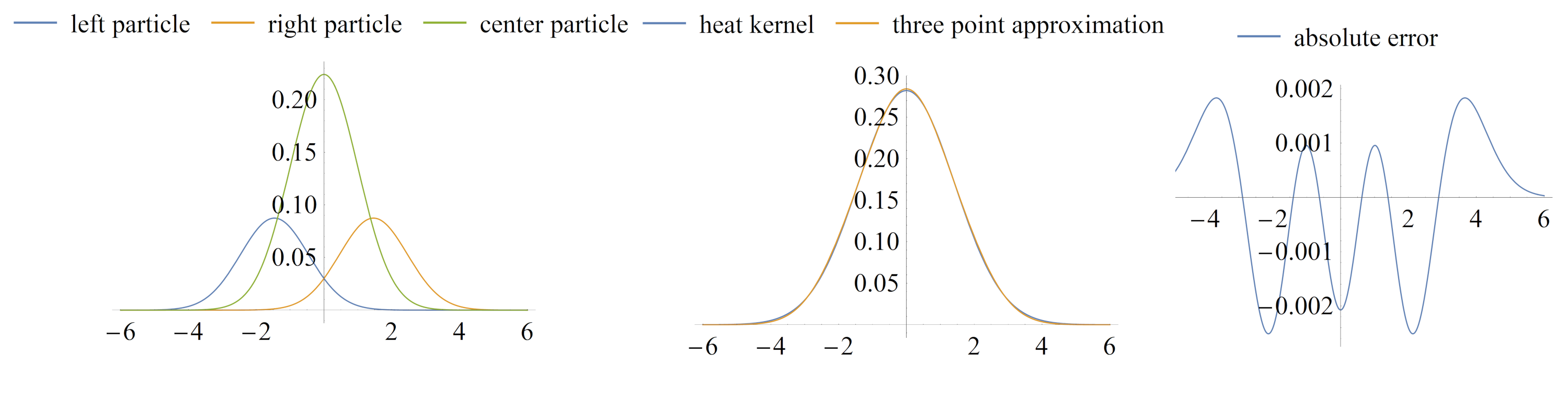}
\caption{Optimal particle split along eigen direction. The left figure plots $3$ components of the split particle. The yellow curve is the new particle in the center, while the green and red curve are the shifted particles. All three have half variance as the density they approximate. The middle figure gives a comparison between a normal distribution (blue) and its approximation (yellow) with $3$ distribution with half variance. The right figure shows the absolute error introduced.}
\label{fig:split_particle}
\end{figure}

We would like to minimize the error introduced measured in the infinity norm, subject to the constraint that the total weight is constant. This reduces this problem to a simple minimization that is independent of the other dimensions of the state space: 
\begin{equation}
\label{eqn:particle_split}
\min_{a,\omega} \max_{x} \left|K_{2d}(x) - ((1-2\omega)K_{d}(x) + \omega K_{d}(x-a) + \omega K_{d}(x+a))\right|
\end{equation}
Where $K_{d}(x)$ is the probability density function of the normal distribution with variance $d$. 
Using numerical optimization, we find $a = 1.03332 \sqrt{d_1}, \omega = 0.21921$. The process is then repeated to other directions, until variance in all directions are sufficiently small.
This process introduces a maximum relative error of $0.73\%$. This is enough to bias a careful analysis of convergence as the size of particles is reduced (and thus, we do not provide convergence results). However, it does offer sufficient accuracy in our simulations since the error introduced has little effect on the linear coupling velocity term. Therefore, the error introduced does not significantly affect the global behavior of the population, as will be shown in the examples below. If higher accuracy is required, one needs to reduce the difference between the variance of the particles. A demonstration of the split is in Fig.~\ref{fig:split_particle}.
\begin{algorithm}[H]
\caption{Particle Split (Eigen-decomposition)}
\label{alg:particle_split}
\begin{algorithmic}
\REQUIRE{Particle with weight $1$ centered at origin with a Diagonal covariance matrix $\bSigma= \text{diag}(\sigma_1 ,\dots,\sigma_d)$}
\STATE {Define the new Diagonal covariance matrix $\bSigma' = \text{diag}(\sigma_1/2,\sigma_2 ,\dots,\sigma_d) $}
\STATE {Find center locations of the off-center particles given by $\bx_c = (1.03332\sqrt{\sigma_1},0,\dots,0)$}
\STATE{Construct center particle centered at origin, with weight $1 - 2\omega$ and covariance $\bSigma'$}
\STATE{Construct left and right particle centered at $\pm\bx_c$, with weight $\omega$ and covariance $\bSigma'$}
\STATE {Apply this method to the new particles until all variance are sufficiently small}
\end{algorithmic}
\end{algorithm}
The above derivation uses the eigendecomposition of the covariance matrix, which could be expensive and undesirable in some applications. While the choice of $a$ and $w$ may not be necessarily optimal for the split, the particle split can be applied to any (non-eigen) direction similarly by substracting appropriate orthogonal projections detailed as follows:

Given a particle centered at origin with weight $1$, and a covariance matrix given by its square root decomposition $\bSigma = \bM \bM^T$, suppose split is needed along direction $\bM_1$, then the children particles are computed as follows: 
\begin{itemize}
\item
Center particle: centered at origin, weight $1-2\omega$, covariance matrix defined by the square root decomposition $\mathbf{N} \mathbf{N}^T$
\item
Left particle: centered at $1.03332 \bM_1$, weight $\omega$, covariance matrix defined by the square root decomposition $\mathbf{N} \mathbf{N}^T$
\item
Right particle: centered at $-1.03332 \bM_1$, weight $\omega$, covariance matrix defined by the square root decomposition $\mathbf{N} \mathbf{N}^T$
\end{itemize}
in which $w=0.21921$, and the matrix $\mathbf{N}$ defined column-wise by: 
\begin{equation}
\label{eqn:split_cov_sqrt}
\mathbf{N}_i = \bM_i - \left(1 - \frac{1}{\sqrt{2}}\right) \frac{<\bM_1,\bM_i>}{<\bM_1,\bM_1>} \bM_1
\end{equation}
in which $<\cdot,\cdot>$ denotes the inner product.
\begin{algorithm}[H]
\caption{Particle Split (square root)}
\label{alg:particle_split_sqrt}
\begin{algorithmic}
\REQUIRE{A particle with weight $1$ centered at origin and a square root decomposition $\bM  \bM^T$ of the Covariance matrix}
\STATE {Find new square root decomposition $\mathbf{N}$ define in (\ref{eqn:split_cov_sqrt})}
\STATE {Find center locations of the off-center particles $\bx_c = 1.03332 \bM_1$}
\STATE{Construct center particle centered at origin, with weight $1-2\omega$ and decomposition $\mathbf{N}$}
\STATE{Construct left and right particle centered at $\pm\bx_c$, with weight $\omega$ and decomposition $\mathbf{N}$}
\STATE {Apply this method to the new particles until all variance are sufficiently small}
\end{algorithmic}
\end{algorithm}
\begin{figure}[H]
\centering
\includegraphics[width =\textwidth]{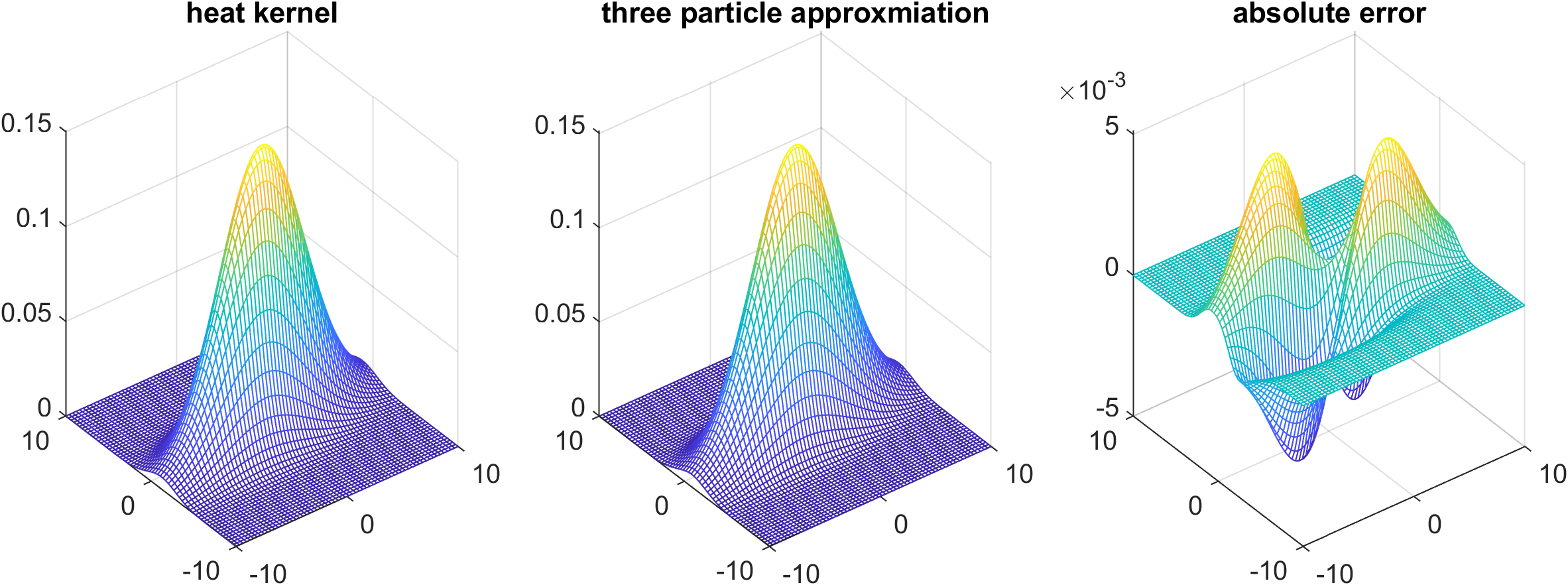}
\caption{Particle split along arbitrary direction.The left figures shows the PDF of a particle with covariance $(16,2;2,3)$. The middle figure and right figure shows the three particle approximation, and absolute error introduced when the particle is split along $(0,1)$ direction, which is not an eigenvector. }
\label{fig:split_particle_sqrt}
\end{figure}
It can be seen in Fig.~\ref{fig:split_particle_sqrt} that while the absolute error introduced is larger than along eigen direction, the split particles still closely approximates the particle they replaced. This split works with the square root representation, and avoids the costly eigendecomposition.
\subsection{Combining particles}
The previous method introduces new particles to the system. Therefore, a method to remove particles is required since otherwise, the number of particles can grow geometrically, which would render this method impossible to use.

The description of the method for combining the particles can be divided into two parts: first, given a set of particles, a method to combine these particles to a single one optimally; second, a strategy to determine which particles should be merged into a single one.

First, consider the method for combining particles.
Except that the total weight does not add to $1$, we are finding a multivariate normal distribution that best approximates the PDF of the old population. Therefore we choose the new particle such that it has the same mean (center location) and covariance matrix of the distribution, and the weight to be the sum of all the particles to be combined, conserving the total weight.

Consider a probability distribution with PDF,
note in the following part: $\boldsymbol{\mu}_i$ is the expectation vector of $i$th particle, whereas $\mu_i$ is the $i$th component of a vector $\boldsymbol{\mu}$:
\begin{equation}
p(\mathbf{x}) = \sum_i \frac{w_i\sqrt{\det(\bSigma_i)}}{(4 \pi)^{d/2}}\exp\left(-\frac{(\mathbf{x} -\boldsymbol{\mu}_i)^T \bSigma_i^{-1} (\mathbf{x} - \boldsymbol{\mu}_i)}{2} \right) .
\end{equation}
Then the expectation of the new distribution $\mathbf{\mu}$ is the weighted average of $\mathbf{\mu}_n$: 
\begin{equation}
\label{eqn:particle_combine_mean}
\boldsymbol{\mu} = \frac{\sum_n w_n \boldsymbol{\mu}_n}{\sum_n w_n}.
\end{equation}
The $ij$th entry of the covariance matrix $\bSigma$ is given by 
\begin{equation}
\label{covarianceijstep1}
\bSigma_{ij} = E[(X_i - \mu_i)(X_j - \mu_j)] = \sum_n \frac{w_n}{\sum_m w_m} E[(Y_{ni} -\mu_i )(Y_{nj} - \mu_j)],
\end{equation}
where $Y_n$ is a multivariate normal random variable subject with expectation $\boldsymbol{\mu}_n$ and covariance matrix $\bSigma_n$. Hence the expectation in \ref{covarianceijstep1} can be evaluated as: 
\begin{align}
E[(Y_{ni} -\mu_i )(Y_{nj} - \mu_j)] &= E[(Y_{ni} - \mu_{ni} + (\mu_{ni} - \mu_i))(Y_{nj} - \mu_{nj} + (\mu_{nj} - \mu_j))]
\\
&=E[(Y_{ni} - \mu_{ni})(Y_{nj}-\mu_{nj})] + (\mu_{ni} - \mu_{i})E(Y_{nj} - \mu_{nj}) \nonumber
\\& + (\mu_{nj} - \mu_{j})E(Y_{ni} - \mu_{ni}) + (\mu_{ni} - \mu_i)(\mu_{nj} - \mu_j) 
\\
&={(\bSigma_{n}})_{ij} + (\mu_{ni} - \mu_i)(\mu_{nj} - \mu_j),
\end{align}
in which $\mu_{nj}$ is the $j$th component of $\boldsymbol{\mu}_n$, and $(\bSigma_{n})_{ij}$ is the $ij$th entry of $\bSigma_n$.
Therefore, we find the covariance matrix $M$ is given by 
\begin{equation}
\bSigma_{ij} =  \frac{1}{\sum_m w_m} \sum_n w_n\left((\bSigma_{n})_{ij} + (\mu_{ni} - \mu_i)(\mu_{nj} - \mu_j)\right).
\end{equation}
In matrix form: 
\begin{equation}
\label{eqn:particle_combine_covariance}
\bSigma = \frac{1}{\sum_m w_m} \sum_n w_n\left(\bSigma_{n} + (\mathbf{\mu}_n - \mathbf{\mu})(\mathbf{\mu}_n - \mathbf{\mu})^T\right)
\end{equation}
 
Then we describe the method to determine which particles should be combined. We choose criteria for combining particles based on their distance being sufficiently close. Since we are more concerned about the computation time for finding combinations rather than have the minimum count of new particles, here we describe a method that scales linearly in the count of particles.

To try combining a particle with others, one can first find all neighbors within a fixed radius, then combine with these neighbors, and iterate through all particles. While apparently finding all particles within a radius threshold needs to compute pairwise distance, it turns out that this fixed-radius near neighbors problem achieves linear scaling in the number of particles \cite{BENTLEY1977}. Unfortunately, there is a lack of implementation for these algorithms for dimensions higher than $3$. Here we describe a much-simplified version that only finds a subset of all neighbors within a fixed radius, since we are only interested in reducing particles quickly, not optimally.

Let the radius threshold for combining particles be $r$. First, we divide the state space into cubic-shaped buckets, with set radius less than $r$. Then, we iterate through particles to find which buckets they are in. After that, we iterate through each nonempty bucket, and combine particles inside to a single particle. Since the particles are only present in a tiny subset of the buckets, the buckets are initialized in a hash table to take advantage of this sparsity, with their bucket index as keys. While the length of the keys grow linearly with the spatial dimension $d$, the time taken to evaluate the hash function can be considered as a constant. Under this assumption, the average case write and access time for a hash table of the occupied bins are practically constant. Hence we achieve practically linear scaling in the number of particles, and constant scaling in the number of dimensions if the number of particles is fixed.

\begin{algorithm}
\caption{Particle Combine}
\label{alg:particle_combine}
\begin{algorithmic}
\STATE {Divide the state space into sufficiently small grids, and find which grid block each particle have its center located (using some Hash function)}
\FOR {each grid block containing more than one particle}
\STATE {combine particles with center location defined as (\ref{eqn:particle_combine_mean}), and covariance matrix defined as (\ref{eqn:particle_combine_covariance})}
\ENDFOR
\end{algorithmic}
\end{algorithm}
\newpage
\section{Simulation results}
\subsection{Motivating example: Van der Pol oscillators}
\begin{figure}[h]
\centering
\includegraphics[width=\textwidth]{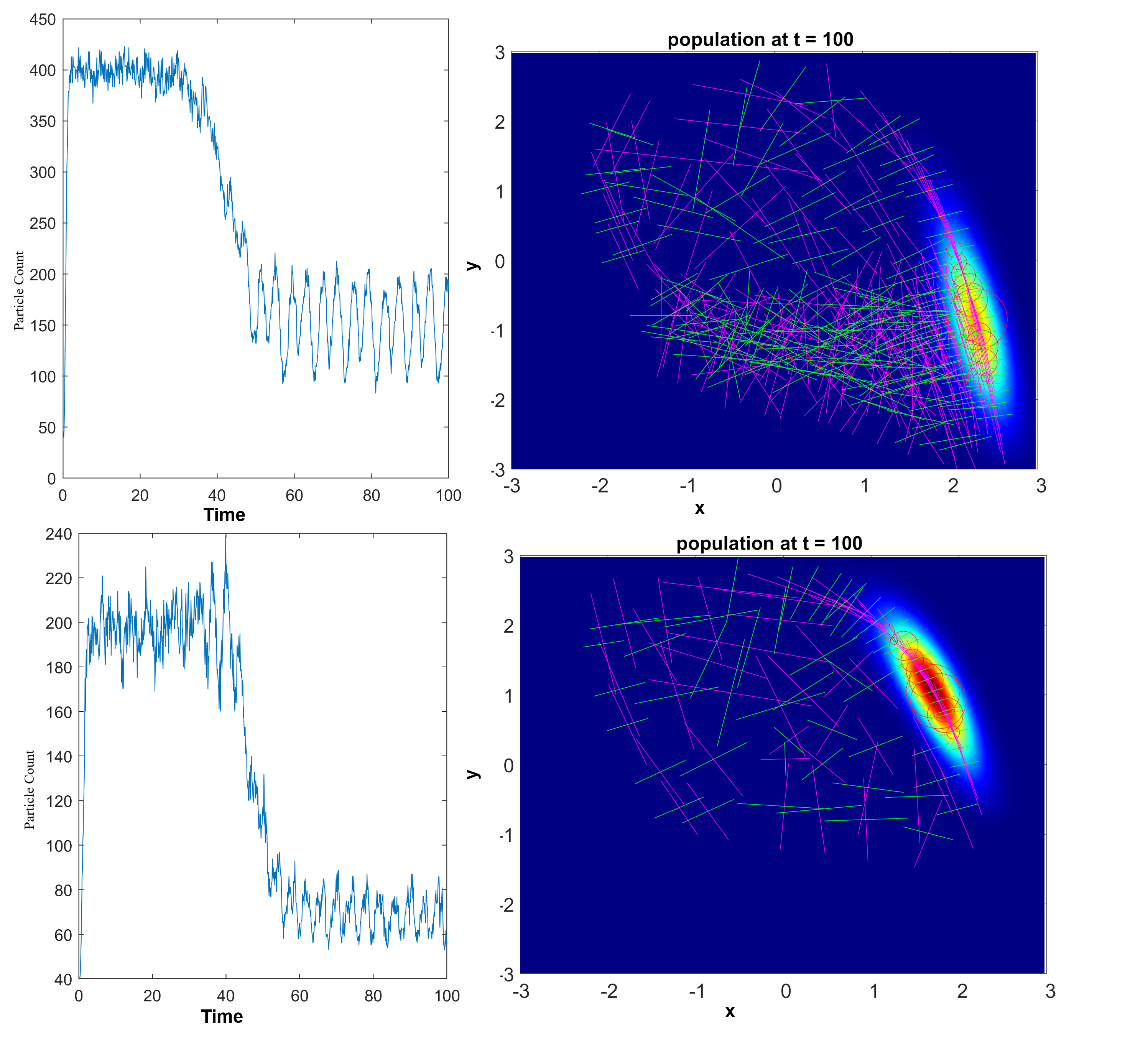}
\caption{Motivating example: The population Van der Pol oscillator with $\mu = 1.5, \alpha=0.5$ and diffusion coefficient $k=0.1$ (top two) and $k=0.05$ (bottom two). In the $2$ graphs to the right: the red circles indicate the weight of individual particles, the magenta, and green lines are the eigenvectors of the covariance matrix with corresponding eigenvalues, while the heat map is the population density at $t=100$. Note that the color scales of the two right graphs are different. The $2$ graphs to the left plot the number of particles versus time. Movies are provided in the supplement.
}
\label{fig:vdp_panels}
\end{figure}

First consider a motivating example of a population of Van der Pol oscillators, as an illustration of the method. Recall (\ref{eqn:fokker_planck}):

\begin{equation}
\frac{\partial u}{\partial t} =  \nabla \cdot \bK \nabla u - \nabla \cdot (\mathbf{v} u),
\end{equation}
where, $\mathbf{v} = \mathbf{v}_d + \mathbf{v}_c$. In the case of the Van der Pol model, the oscillator dynamic velocity field $\mathbf{v}_d = \partial \mathbf{x}/\partial t$ is defined by: 
\begin{equation}
\begin{bmatrix}
v_1\\v_2
\end{bmatrix}=
\begin{bmatrix}
\mu (x_1 - \frac{1}{3}x_1^3 - x_2)\\
\frac{1}{\mu} x_1
\end{bmatrix}.
\end{equation}
The coupling velocity $\mathbf{v}_c = \mathbf{v}_c(u)$ is defined by
\begin{equation}
\mathbf{v}_c(u) =
\alpha
\int_{\mathbb{R}^2} y_{1} u(\mathbf{y}) d \mathbf{y},
\end{equation}
in which $\alpha$ is a coupling coefficient. This is commonly referred as a \emph{mean-field} coupling. We set up the initial population along the limit cycle of this oscillator without coupling. 
The results from a sample simulation is shown in Fig.~\ref{fig:vdp_panels} 
, which serves as an illustrative example on how the APPD represents the population and noise over the process of simulation.

We choose the Van der Pol oscillator to serve as an illustrative example for the usage of our algorithm, as the state space is $2$-dimension and can be plotted without projection. In the example used in Fig.~\ref{fig:vdp_panels}, we select an uncoupled initial condition along the limit cycle and observe the effects of coupling. It can be seen that both populations are coupled; however, there is a slight phase difference due to the noise level difference. It started at a high value since the population is set to start in an uncoupled state. As the coupling synchronizes the population, the population density is concentrated in a smaller region, and the number of particles drops accordingly. Since the computational cost is proportional to particle count, our method is more efficient if there is stronger coupling and weaker diffusion.

Additionally, we note that synchronization occurs unevenly and that once enough particles start to synchronize, many of the other particles quickly follow suit. This occurs in a similar way to the phase transitions seen in coupled oscillator theory~\cite{forger2017biological}. 
As the population synchronizes, the number of particles shrinks. 

\subsection{Hodgkin-Huxley model with threshold coupling}
The Hodgkin-Huxley model is a model that describes the electrical activity of a neuron based on ionic currents. Coupled neurons form the basis of computational neuroscience.
These models are highly nonlinear and have multiple attractors, which creates challenges for simulation. Additionally, the model exhibits complex behaviors that are highly dependent on the noise level, such as noise-induced synchronization~\cite{wang2000coherence,BRUNEL2006}, and a noise-induced coexistence of firing and resting neurons~\cite{bashkirtseva2015stochastic}. These interesting macroscopic phenomena require an accurate simulation of noise to reproduce. We use the Hodgkin-Huxley model in two ways. First, we show that the APPD accurately reproduces the mentioned phenomena, which can be further separated into two cases. Then, we compare it against direct Monte Carlo simulation to show that the method is accurate and fast.

The Hodgkin-Huxley model used here is described as follows:
We take $\mathbf{x} = (0.01V,m,n,h)$, where $V$ is the membrane potential in millivolts, $m,h$ are proportion(for each cell) of activating and inactivating subunits of the sodium channel, and $n$ corresponds to that of potassium channel subunits. Here $V$ is scaled by $0.01$ such that all the values have the same order of magnitude, where our method is most efficient.

The model equations are as follows:
\begin{align}
\label{eqn:hh_conduction}
C\frac{\partial}{\partial t}V &= -G_\text{Na}m^3h(V-V_\text{Na})-G_\text{K}n^4(V-V_\text{K}) \\\nonumber
&- G_L(V-V_L) - G_\text{coupling} (V - V_\text{coupling})- I_\text{app}\\
\frac{\partial m}{\partial t} &= \alpha_m(V)(1-m) - \beta_m(V)m\\
\frac{\partial h}{\partial t} &= \alpha_h(V)(1-h) - \beta_h(V)h\\
\frac{\partial n}{\partial t} &= \alpha_n(V)(1-n) - \beta_n(V)n.
\end{align}
With the following equations for the subunit dynamics:
\begin{align}
\alpha_m(V) &= 0.1 \frac{V-25}{1-\exp(-\frac{V-25}{10})}\\
\beta_m(V) &= 4\exp(-\frac{V}{18})\\
\alpha_h(V) &= 0.07\exp(-\frac{V}{20})\\
\beta_h (V) &= \frac{1}{1+\exp(-\frac{V-30}{10})}\\
\alpha_n(V) &= 0.01\frac{V-10}{1-\exp(-\frac{V-10}{10})}\\
\beta_n(V) &= 0.125\exp(-\frac{V}{80}).
\end{align}
in which $G_\text{coupling}$ is conductance due to coupling, as explained below. 

In a neuron population, it is usually assumed that individual neurons are independent, except when a neuron firing occurs, and post-synaptic neurons are coupled with the firing neuron. In a population density model, the neurons are indistinguishable except for their state, therefore an all-to-all coupling (or probabilistic coupling) between neurons is assumed. However, contrary to a neural mass model where coupling strength is determined from the average membrane potential~\cite{nykamp2000population,david2003neural}, the more realistic threshold coupling can be implemented. Specifically, $G_\text{coupling}$ is proportional to the flow rate across the hyperplane $V = V_\text{threshold}$ from $V < V_\text{threshold}$ to $V > V_\text{threshold}$ side. Since the threshold value is chosen where the transition is fast, the contribution of diffusion to coupling can be ignored, and the flow rate is given by the velocity of particles multiplied by the marginal density along $V$. In the actual implementation, we used the averaged flow rate over the previous timestep as the flow rate for the current step, which, considering that actual coupling is not instantaneous, is still a reasonable assumption. The threshold firing potential is chosen at $V_\text{threshold} = 45 \si{\milli\volt}$.

Since the subunit variables are proportions of some quantity, their domain is restricted in $[0,1]$. The dynamics have inward-pointing velocity on the boundary. Therefore the center of particles will not leave the domain without splitting. However, as the particles have volume, the support of centers can be outside this domain due to the effects of Gaussian noise. \cite{alos2002stochastic} and \cite{krepysheva2006space} discuss more detailed handling of the noise near boundary conditions. Nonetheless, we find that the following alternative method can still produce sufficiently accurate results.

The noise can result in off-domain points in two ways: When a split occurs and when checking an off-center point of a particle in the LSKF method. If a split results in the center of a particle outside the domain, the center is shifted to the closest location inside the domain. If an off-center point is outside the domain while the center is inside, we choose the offset of the point from the particle's center to be its negative, which will be inside the domain and on the level set.

The value of the constants are listed as follows: $C=1\si{\micro\farad\per\square\centi\metre}, G_\text{Na} = 120\si{\micro\ampere\per\milli\volt\square\centi\metre}, E_\text{Na} = 115\si{\milli\volt}, G_\text{K} = 36\si{\micro\ampere\per\milli\volt\square\centi\metre}, E_\text{K} = -12\si{\milli\volt}, G_L = 0.3\si{\micro\ampere\per\milli\volt\square\centi\metre}, E_L=10.613\si{\milli\volt}, I_\text{app} = 10 \si{\micro\ampere}$. The coupling potential $V_\text{coupling}=-35 \si{\milli\volt}$ for the inhibitory case, and $V_\text{coupling}=50\si{\milli\volt}$ for the excitatory case, keeping in mind that in this original Hodgkin-Huxley model, the neuron rests near $0  \si{\milli\volt}$. We choose the diffusion to be homogeneous after $V$ is scaled, with a diffusion matrix $\mathbf{K} = k \mathbf{I}_d$.

The coupling is defined analogously to neuron firing: that is, when a neuron reaches a threshold membrane potential from below, the neuron sends a signal to all post-synaptic neurons. Consequently, $G_\text{coupling}$ in (\ref{eqn:hh_conduction}) is defined by $G_\text{coupling} = 20Qc$, where constant $c$ is the coupling coefficient, $Q$ is the flow rate of neurons as a proportion of total population per millisecond. (Therefore, $Q$ has the unit $\si{\per\milli\second}$.)

\begin{figure}[H]
\includegraphics[width =\textwidth]{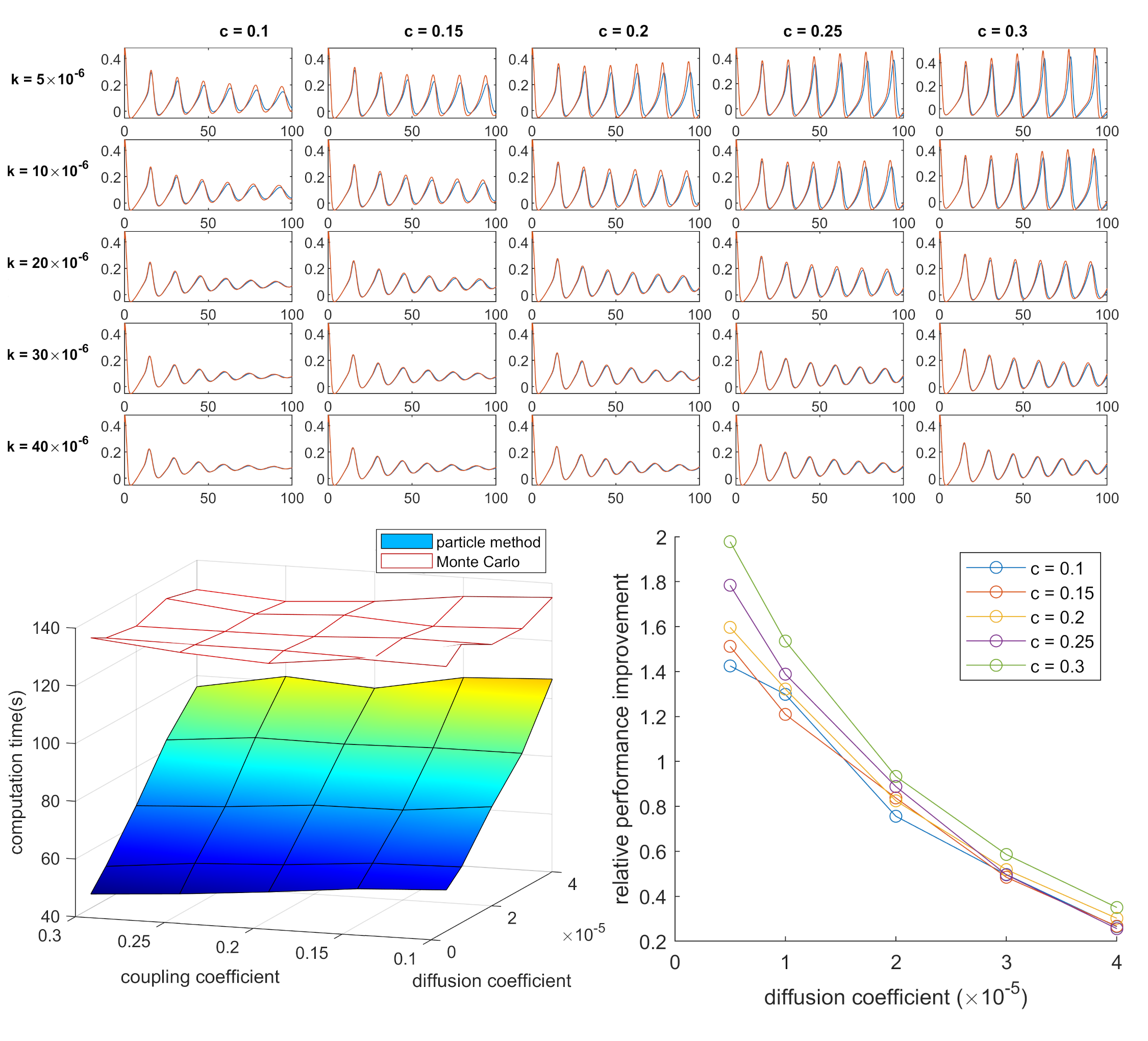}
\caption{Testing the accuracy and efficiency of the method. This figure compares the APPD (blue) result against a Monte Carlo (orange) method with $41080$ neurons for excitatory neurons over a range of parameters. The top panel plots the Average membrane potential $V$(left $5$ columns) versus time in microseconds. There are $25$ different instances with varying levels of diffusion and coupling strength: From left to right, the coupling coefficient $c$ increases from $0.1$ to $0.3$ with a stepsize of $0.05$. From top to bottom, the diffusion coefficient increases, with $k=0.5,1,2,3,4 \times10^{-5}$ respectively. The lower panel plots the computation time taken for the computation: the left panel shows the computation time, whereas the right panel shows the relative performance improvement over the Monte Carlo method}
\label{fig:HH_panels_B}
\end{figure}
We compare our simulations to direct Monte Carlo simulations. Note we choose the number of elements $N = 41080$ for the Monte Carlo simulation, and the improvement in performance would be more significant if a larger number of elements is used. For the range of parameters chosen and our implementation of the Monte Carlo method, the APPD is faster than the direct Monte Carlo method, as shown in Fig.~\ref{fig:HH_panels_B}. It is important to note that while the Monte Carlo method takes about a constant amount of time for all simulation cases, there are large variations for the APPD. At a high noise, low coupling scenario ($k=4\times 10^{-5},c=0.1$) where the population is asynchronous at the end of the simulation, the APPD is faster than the Monte Carlo simulation by $20\%$. Whereas for a low noise strong coupling scenario ($k=0.5\times10^{-5},c=0.3$), the population remains synchronized, and the APPD is $200\%$ faster, as shown in the lower-right panel of Fig.~\ref{fig:HH_panels_B}. This again shows that the APPD is most suitable to simulate neurons with strong coupling at a lower noise level and that the APPD can adapt to create computational efficiencies.
\section{Conclusions}
The APPD method provides a fast and accurate method to study population-level behaviors while accounting for noise and coupling in biologically realistic models. We compared our method against direct Monte Carlo simulation and tested its ability to reproduce complex macroscopic phenomena across a range of parameters.

While the APPD is conceptually inspired by the particle method by Stinchcombe and Forger~\cite{STINCHCOMBE2016}, the method presented here has two significant improvements: 1)~particles are asymmetric, allowing it to better track the population; 2)~computation of the noise term is computed by modifying a single particle instead of interactions between nearby particles, which reduces computation, and allows almost perfect parallel implementation.
When compared against existing asymmetric particle methods~\cite{rossi1996resurrecting,rossi2005achieving,xie2014adaptive,berchet2021adaptive}, our method is the first to be applied to a $4$-dimensional problem instead of the typically considered $2$-dimensional problems. Such an extension is nontrivial as it depends on the improved numerical stability of our method by tracking a \emph{square-root} instead of the covariance. 

Looking at the examples, we checked whether we achieved the goal to develop a middle ground between the direct simulation and existing population density-based methods. From the form of (\ref{eqn:level_set_ODE_num_matform}), it is evident that the dynamics of the oscillator are retained, and no manual dimension reduction is required. The choice of the kernel (\ref{eqn:kernel_function}) enabled an accurate representation of the effect of noise. Through the examples, we implemented both all-to-all-coupled oscillators and threshold-coupled oscillators; both resemble the direct simulation.

The advantages of the APPD are the versatility through different models and a range of parameters without any model-specific simplifications, and the robustness that the macroscopic behavior matches the direct simulation (Fig.~\ref{fig:HH_panels_B}). These observations lead us to conclude that the APPD can be applied to analyze the complex macroscopic behaviors of a population of coupled noisy oscillators, which is in contrast to many population density approaches that require a model-specific understanding of its behavior to be implemented.

\section{Acknowledgment}
We would like to thank Adam Stinchcombe for his input through the development of the method and shaping of the manuscript. This project is partially supported by the NSF1714094 grant and HFSP program grant RGP 0019-2018.

\bibliography{ref}
\bibliographystyle{elsarticle-num}
\end{document}